\documentclass[letterpaper, 10 pt, conference]{ieeeconf}  

\IEEEoverridecommandlockouts                              
\usepackage[utf8]{inputenc}

\usepackage{graphics}
\usepackage{graphicx,epstopdf}
\usepackage{amssymb,amsmath,bm,dsfont}
\usepackage{ams fonts}
\usepackage[rm]{subfigure}
\usepackage{xcolor}

\newtheorem{theorem}{Theorem}

\DeclareMathOperator*{\argmin}{argmin}

\DeclareMathOperator{\interior}{int}

\newcommand{\norm}[1]{\ensuremath{\left\| #1 \right\|}}

\newcommand{\calC}{\ensuremath{\mathcal{C}}}

\newcommand{\calF}{\ensuremath{\mathcal{F}}}

\newcommand{\calK}{\ensuremath{\mathcal{K}}}

\newcommand{\calN}{\ensuremath{\mathcal{N}}}
\newcommand{\calO}{\ensuremath{\mathcal{O}}}

\newcommand{\calX}{\ensuremath{\mathcal{X}}}

\newcommand{\bI}{\ensuremath{\bm{I}}}
\newcommand{\bJ}{\ensuremath{\bm{J}}}

\newcommand{\bp}{\ensuremath{\bm{p}}}

\newcommand{\bz}{\ensuremath{\bm{z}}}

\newcommand{\setR}{\ensuremath{\mathbb{R}}}

\def\st/{\textsuperscript{st}}
\def\nd/{\textsuperscript{nd}}
\def\rd/{\textsuperscript{rd}}
\def\th/{\textsuperscript{th}}

\newcommand{\del}{\ensuremath{\partial}}
\newcommand{\ones}{\ensuremath{\bm{\mathds{1}}}}

\title{\LARGE \bf
Source Seeking in Unknown Environments with Convex Obstacles
}

\author{Bruno A. Ang\'elico$^{1}$, Luiz F. O. Chamon$^{2}$, Santiago Paternain$^{2}$, Alejandro Ribeiro$^{2}$, and George J. Pappas$^{2}$%
\thanks{$^{1}$Dept. of Telecomm. and Control Engineering, Escola Polit\'ecnica, University of S\~ao Paulo, Brazil~({\tt\small angelico@lac.usp.br}). B.A.~Ang\'elico is supported by Funda\c{c}\~ao de Amparo \`a Pesquisa do Estado de S\~ao Paulo (FAPESP), grant 2018/12359-7.}%
\thanks{$^{2}$Dept. of Electrical and Systems Engineering, University of Pennsylvania, USA \big(\tt\small \{luizf, spater, aribeiro, pappasg\}@seas.upenn.edu\big)}%
}

\begin{document}

\maketitle
\thispagestyle{empty}
\pagestyle{empty}

\begin{abstract}

Navigation tasks often cannot be defined in terms of a target, either because global position information is unavailable or unreliable or because target location is not explicitly known \emph{a priori}. This task is then often defined indirectly as a source seeking problem in which the autonomous agent navigates so as to minimize the convex potential induced by a source while avoiding obstacles. This work addresses this problem when only scalar measurements of the potential are available, i.e., without gradient information. To do so, it construct an artificial potential over which an exact gradient dynamics would generate a collision-free trajectory to the target in a world with convex obstacles. Then, leveraging extremum seeking control loops, it minimizes this artificial potential to navigate smoothly to the source location. We prove that the proposed solution not only finds the source, but does so while avoiding any obstacle. Numerical results with velocity-actuated particles, simulations with an omni-directional robot in ROS+Gazebo, and a robot-in-the-loop experiment are used to illustrate the performance of this approach.

\end{abstract}


\section{Introduction}\label{sec:sec_intro}

Navigation tasks are often defined by specifying a target position that the autonomous agent must attain~~\cite{Choset, LaValle}. This, however, requires access to global information that is scarce in many applications. For instance, when navigating underwater or in-door, global position information is either not available, due to the lack of GPS signal, or unreliable, due to the high cost of obtaining long-term accurate inertial measurements. What is more, the target location need not be known explicitly \emph{a priori}, e.g., when locating the source of an oil spill~\cite{Senga} or radio source localization~\cite{Twigg}. A common way solution for these scenarios is to define the robot's task indirectly as a source seeking problem.

Source seeking can be abstracted as the problem of minimizing~(maximizing) a convex~(concave) potential induced by the source. Though the potential is unknown, its values and more importantly its gradient can be estimated from sensor measurements. A gradient controller can then be used to drive the robot to the potential's minimum~(maximum) as long as it navigates an open convex space, i.e., an environment without obstacles~\cite{Zhang2007, Ghods2007}. This solution is appealing not only because it does not require prior explicit knowledge of the goal, but because of its simplicity and the fact that they require only local information.

Still, if this method is to be useful is real-world applications, two issues must be addressed: (i)~how to obtain good enough gradient estimates from local, scalar measurements and (ii)~how to navigate cluttered environments. The main approaches used to tackle issue~(i) involve variations of the classical Robbins-Monro algorithm, such as random directions for stochastic approximation~(RDSA)~\cite{Azuma2012, Atanasov2012, Ramirez-Llanos2019} or sinusoidal extremum seeking control~(ESC)~\cite{Zhang2007, Ghods2007, Cochran2009, Ghods2010}. Albeit effective, the stochastic nature of the former may lead to non-smooth trajectories as the robot takes both random steps to sample values of the potential and steps in the direction of the estimated gradient. Although sinusoidal ESC does not suffer from this issue, the fact remains that neither methods are guaranteed to work in the presence of obstacles.

Artificial potentials can be used to navigate cluttered environments using gradient dynamics~[issue~(ii)]. The idea is to combine the attractive potential defining the agent's goal with repulsive potentials representing the obstacles. Yet, careful combination is required to avoid creating spurious local minima that would prevent the robot to complete the original task. Nevertheless, by following the gradient of Koditschek-Rimon potentials or navigation functions, an agent is guaranteed to the navigate to the minimum of a strongly convex potential while avoiding strongly convex, non-intersecting obstacles~\cite{koditschek1990robot, Rimon1992, Filippidis, Paternain2018}.

This paper tackles the practical problem of source seeking from scalar measurements in a cluttered environment. To do so, we use local measurements from the objective potential and~(partial) knowledge of obstacles to construct an artificial potential~(Section~\ref{sec:sec_NF}) that the agent navigates using an ESC loop~(Section~\ref{sec:sec_ESC}). Note that neither the target nor the environment needs to be known \emph{a priori}: the artificial potential is evaluated online from sensor measurements. We prove that under mild conditions this autonomous system is guaranteed to attain its goal while avoiding collisions~(Section~\ref{sec:sssm}) and illustrate its performance for an omni-directional robot navigating around spherical obstacles~(Section~\ref{sec:sec_results}).

\section{Problem Formulation}

As we mentioned in the introduction, we pose the problem of source seeking in a cluttered environment as that of minimizing an unknown scalar potential through a trajectory defined in a non-convex space. Explicitly, let~$\calX \subset \setR^2$ be a non-empty, compact, convex set delineating an environment that contains non-intersecting~(possibly unknown) obstacles represented by the open convex sets~$\calO_i \subset \calX$ with non-empty interior and smooth boundaries~$\del\calO_i$. The \emph{free space}, i.e., the space of points that the agent can occupy, is therefore defined as
\begin{equation}\label{E:free_space}
	\calF = \calX \setminus \bigcup_i \calO_i.
\end{equation}

Our goal is to reach a source located at~$\bp^\star \in \interior(\calF)$~(unknown) that induces a strongly convex potential~$f_0: \calX \to \setR_+$ such that~$\bp^\star = \argmin_{\bp \in \calX} f_0(\bp)$. Observe that many physical phenomena, such as those involving spherical wave propagation or~$1/r^2$ decay, induce strongly convex potentials. We wish to do so using only scalar measurements of the potential~$f_0$ and while remaining in the free space~$\calF$ at all times. Formally, our goal is to generate a trajectory~$\bp(t)$ for the agent such that~$\bp(t) \in \calF$, for all~$t \in [0,\infty)$, $\lim_{t \to \infty} \bp(t) = \bp^\star$, and~$\dot{\bp}(t)$ is a function only of~$f_0(\bp(\tau))$ for~$\tau \leq t$, i.e., without gradient measurements.

Notice that the space~$\calF$ in which the robot will navigate is punctured and therefore non-convex. Hence, even if the gradient of~$f_0$ were known, simple gradient dynamics is not guaranteed to stay within~$\calF$. This problem can be circumvented by constructing an artificial potential such as the Koditschek-Rimon navigation function from local measurements~\cite{koditschek1990robot, Paternain2018}. However, the question remains of whether gradient estimate errors would lead to collisions. Though this issue was addressed for stochastic estimates in~\cite{paternain2016stochastic}, stochastic approximation methods may lead to highly discontinuous trajectories.

To address these issues, we start by describing the two building blocks of our solution: navigation functions and sinusoidal ESC~(Section~\ref{sec:prelim}). We then introduce our solution and prove that, under mild conditions, any particle following this policy will approach the source without colliding with obstacles~(Section~\ref{sec:sssm}). Finally, we illustrate the performance of our method using an omni-directional robot navigating a space with spherical obstacles~(Section~\ref{sec:sec_results}).

\section{Technical Background}\label{sec:prelim}

In this section, we briefly introduce the two building blocks we use to address the problem of source seeking in cluttered environments: navigation functions and sinusoidal ESC. More detailed discussions of these topics can be found, e.g., in~\cite{Paternain2018} and~\cite{Krstic2003} respectively.

\subsection{Navigation Functions}\label{sec:sec_NF}

Navigation functions are artificial potentials that enable collision-free navigation of spaces with obstacles. Explicitly, a map~$\nu: \calF \to [0,1]$ is said to be a \emph{navigation function} towards~$\bp^\star$ in~$\calF$ if (i)~$\nu \in \calC^2$~(twice continuously differentiable), (ii)~$\nu$ is polar at~$\bp^\star$, (iii)~$\nu$ is Morse, and (iv)~$\del \calF = \nu^{-1}(1)$~\cite{koditschek1990robot}. This map is called a navigation function because the dynamical system
\begin{equation}\label{E:dynamics}
	\dot{\bp} = -\nabla \nu(\bp)
\end{equation}
navigates the free space~$\calF$ until it reaches~$\bp^\star$~\cite{Koditschek_Lyap}.

In this work, we use the Rimon-Koditschek navigation function~$\varphi$ from~\cite{koditschek1990robot}. This potential is constructed by describing each obstacle~$\calO_i$ as the null sublevel set of a convex function~$\beta_i: \calX \to \setR$. Explicitly, $\calO_i = \{\bp \in \calX \mid \beta_i(\bp) < 0\}$. Such a function always exists since every convex set is the sublevel set of a convex function~\cite{Boyd}. Then, for~$\beta(\bp) = \prod_{i} \beta_i(\bp)$, the Rimon-Koditschek potential is defined as
\begin{equation}\label{eq:eq_NV}
	\varphi(\bp) \triangleq \frac{f_0(\bp)}{\left(f_0^k (\bp) + \beta(\bp) \right)^{1/k}}
		\text{,}
\end{equation}
where~$k > 0$ is a fixed order parameter and~$f_0$ is the potential induced by the source. Note that if the obstacles~($\beta_i$) are known, then~$\varphi(\bp)$ can be evaluated by measuring the potential~$f_0(\bp)$. If they are not, $\varphi$ can be built online as the agent encounters the obstacles. In this case, if the obstacles are ellipsoids, the functions~$\beta_i$ can be estimated by measuring their curvatures~\cite{Paternain2018}.

It has been established that for sufficiently curved worlds, the artificial potential~\eqref{eq:eq_NV} is a navigation function for large enough~$k$~\cite{filippidis2012navigation, Paternain2018}. What is more, even when the gradient in~\eqref{E:dynamics} is replaced by a stochastic approximation, the agent is guaranteed to navigate the free space to~$\bp^\star$ under mild conditions~\cite{paternain2016stochastic}. However, in order to deal with practical robotic applications, we wish to obtain smoother trajectories than the ones produced by stochastic gradient approximations. To do so, we rely on extremum seeking dynamics.

\subsection{Sinusoidal Extremum Seeking Control}\label{sec:sec_ESC}

\begin{figure}[t]
	\centering
	\includegraphics[width=0.6\linewidth]{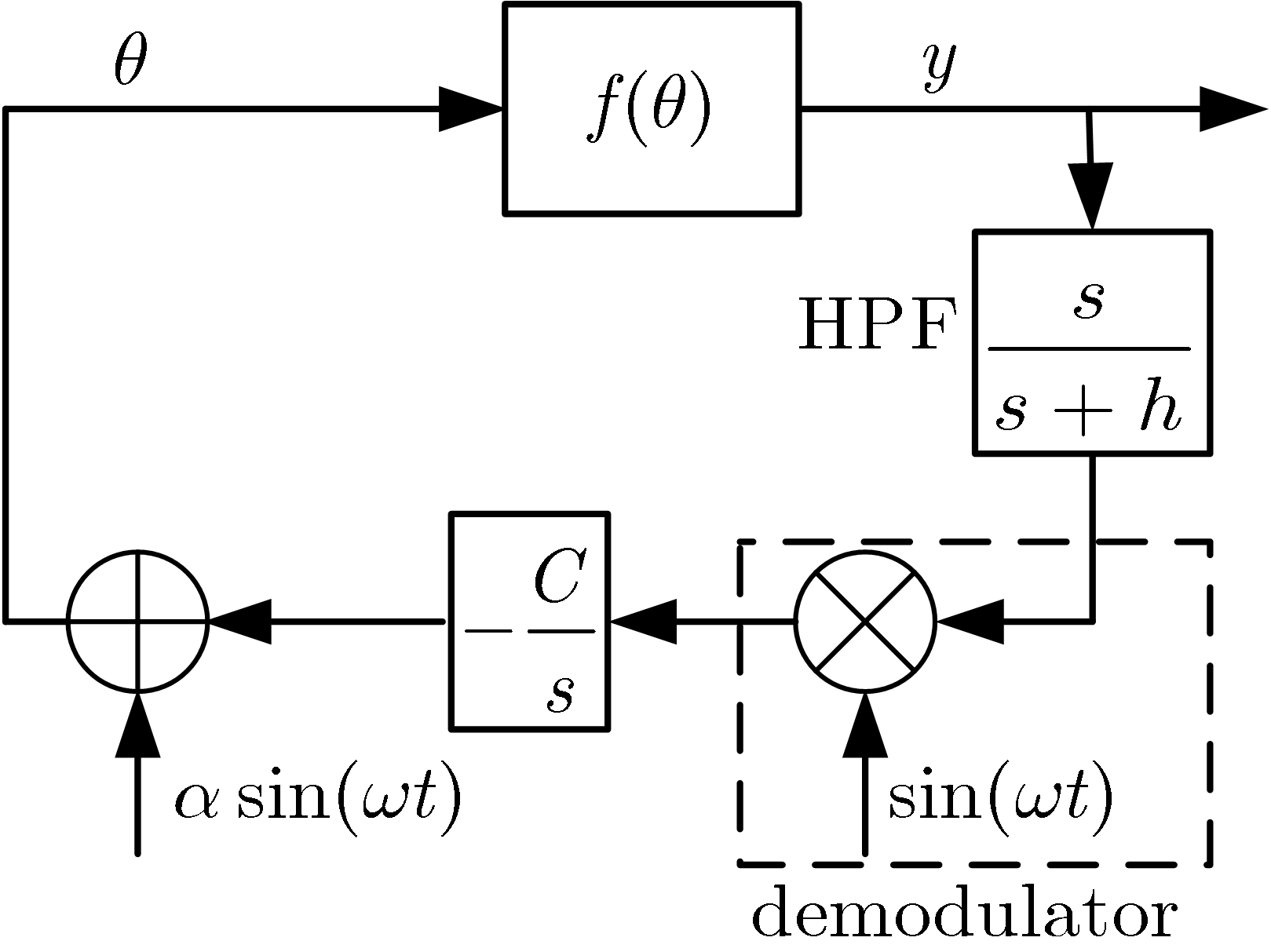}
	\caption{Sinusoidal extremum seeking controller for static map and minimization purposes.}	
	\label{fig:ESC_min_02}
\end{figure}

The idea of ESC emerged in the context of gradient-free optimization~\cite{Leblanc1922}, having since spread to applications in model-free control, adaptive control, and source seeking~\cite{Astrom1995, Tan2010, Scheinker2017,Zhang2007, Ghods2007, Cochran2009, Ghods2010, Raisch2017}. The basic ESC loop is shown in Fig.~\ref{fig:ESC_min_02}. Intuitively, ESC uses a small periodic perturbation~$\alpha \sin(\omega t)$ to explore the landscape of the function~$f$ around~$\theta$. Assuming~$f$ is a smooth convex function, this perturbation creates a periodic response in~$y$ that is in phase with respect to~$\sin(\omega t)$ if~$\theta > \theta^\star$ or out of phase if~$\theta < \theta^\star$, for~$\theta^\star = \argmin_\theta f(\theta)$. The feedback loop then corrects~$\theta$ until it stabilizes at~$\theta^\star$. Another way to understand ESC is by noticing that the high-pass filter~(HPF) and the demodulator are effectively estimating the gradient of~$f$~\cite{Krstic2003}. Naturally, the HPF cut-off frequency must be lower than that of the sinusoidal perturbation~($h < \omega$).

The loop in Fig.~\ref{fig:ESC_min_02} is guaranteed to converge to the minimum~$\theta^\star$ for reasonable choices of parameters and under mild conditions on~$f$~\cite{Krstic2000}. Since, these results have been extended to vector-valued parameters and loops including nonlinear plants dynamic~\cite{Krstic2003}. It is not however straightforward that using ESC to approximate the navigation function dynamics in~\eqref{E:dynamics} will lead to a safe trajectory. Indeed, though this combination will eventually reach the source location~$\bp^\star$, it need not do so through a trajectory that remains in the free space. To guarantee this is the case, we need to show that the ESC loop provides a sufficiently accurate estimate of the navigation function gradient for all~$t$. In the following section, we describe how the navigation function in~\eqref{eq:eq_NV} can be evaluated using local measurement. We then construct an ESC loop to follow its gradient and prove that it in fact does so, showing that this approach is able to navigate a cluttered environment without collision to reach the source of interest.

\section{Source Seeking with Scalar Measurements}\label{sec:sssm}

In this section, we describe the proposed solution for the problem of source seeking with scalar measurements. To do so, we consider a velocity-actuated point mass as in~\cite{Zhang2007}. A straightforward extension to omni-directional robots is explored in the experiments of Section~\ref{sec:sec_results}. Though the presence of more involved dynamics has been considered in the context of navigation functions~\cite{Koditschek_Mech} and ESC~\cite{Krstic2003}, we defer this issue for future work.

\begin{figure}[t]
	\centering
	\includegraphics[width=0.65\linewidth]{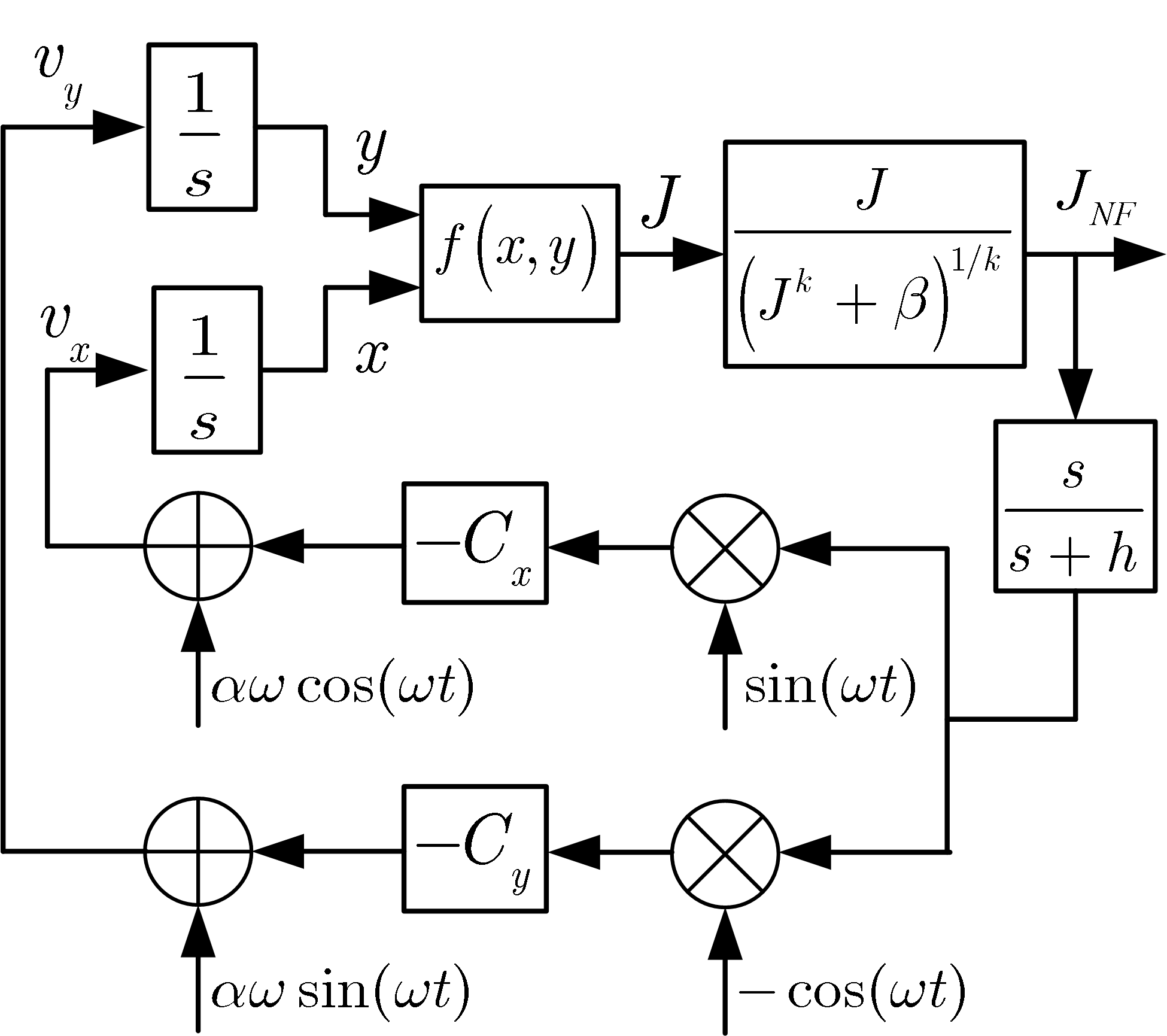}
	\caption{Source seeking with scalar measurement in cluttered environment for velocity-actuated particle whose position is given by~$\bp = (x,y)^T$.}
	\label{fig:ES_NV}
\end{figure}

A diagram of the hybrid navigation function/ESC solution is presented in Fig.~\ref{fig:ES_NV}. The agent senses the source by measuring its potential at the current location~$\bp = (x,y)^T$. It then evaluates the navigation function~$\varphi$ in~\eqref{eq:eq_NV} by incorporating partial knowledge of the obstacles, i.e., using~$\tilde{\beta}(\bp) = \prod_{i \in \calK} \beta_i(\bp)$ where~$\calK$ is the set obstacles that the agent has encountered so far. Proceeding, this value is fed to a two-dimensional version of the ESC loop from Fig.~\ref{fig:ESC_min_02}, constructed as in~\cite{Zhang2007} by using orthogonal perturbations, namely~$\sin(\omega t)$ and~$\cos(\omega t)$. Finally, the gradient estimate from the ESC loop is feedback to the agent actuators.

The ESC loop in Fig.~\ref{fig:ES_NV} is guaranteed to converge to a small neighborhood around~$\bp^\star$. Indeed, \cite{Zhang2007} showed that the error~$\varphi(\bp(t)) - \varphi(\bp^\star)$ converges to a~$\calO(\alpha^2 + 1/\omega^{2})$-neighborhood of zero for reasonable choices of the gains~$C_x,C_y$. Since~$\varphi$ is polar, this implies that~$\lim_{t \to \infty} \bp(t) \in \calN^\star$, a small neighborhood of~$\bp^\star$. Moreover, the approximation error can be reduced by using small, high frequency perturbations. This, however, does not guarantee that~$\bp(t) \in \calF$ for all~$t \geq 0$, i.e., that the agent does not collide with obstacles. To show this is the case, we must demonstrate that ESC not only converges to the minimum of~$\varphi$, but does so by closely following its gradient. In other words, demonstrate that ESC in fact approximates the dynamics in~\eqref{E:dynamics}. The following theorem characterizes the trajectory~$\bp(t)$ and shows that under mild conditions, the agent approaches the source while avoid obstacles.

\begin{theorem} \label{theo:theo_Luiz}
Consider the dynamical system in Fig.~\ref{fig:ES_NV} and assume that~$C_x = C_y = C \ll h \ll \omega$ and~$\alpha \ll 1$, i.e., that the loop gains and the amplitude of the perturbation are small whereas the HPF cut-off and perturbation frequencies are large. Then, there exists a safety guard~$\sigma = \calO(\alpha C/\omega)$ such that if the agent considers virtual obstacles inflated by~$\sigma$, it will navigate collision-free to a neighborhood of~$\bp^\star$.
\end{theorem}

\begin{proof}
See appendix.
\end{proof}

\section{Performance Results}\label{sec:sec_results}

Three different scenarios were considered. In the first one, a Matlab simulation with a simple velocity-actuated particle model was performed. In the second case, a more realistic simulation with a commercial Quadruple Mecanum Robot was carried out using ROS+Gazebo. Lastly, an experiment considering a real dynamics with a robot-in-the-loop was considered. In all cases, the obstacles were inflated approximately by the radius of the robot, and the source potential was modeled as:-
\begin{equation}
	J = q_x (x-x_s)^2 + q_y (y - y_s)^2
    \label{eq:eq_J}
\end{equation}

\subsection{Velocity-actuated Particle}

The parameters of the ESC were chosen such that $\omega = 40$ rad/s, $\alpha = 0.07$, $C_x = C_y = -10$ and $h=20$. The world consists in a circled environment with radius 3 m and five round obstacles with radius 0,25 m at the points $(-1.0,0.0)$, $(-0.2,1.2)$, $(1.0,0.7)$, $(1.0,-1.0)$, $(-0.5,-1.0)$ , respectively. For the source potential, it was considered $q_x =  q_y = 1.0$, and for the navigation function, it was set $k = 6$. The starting point of the robot was se to $(0, 2.5)$. 

\begin{figure}[hbt!]
	\centering
	\subfigure[][]{\includegraphics[width=.7\columnwidth]{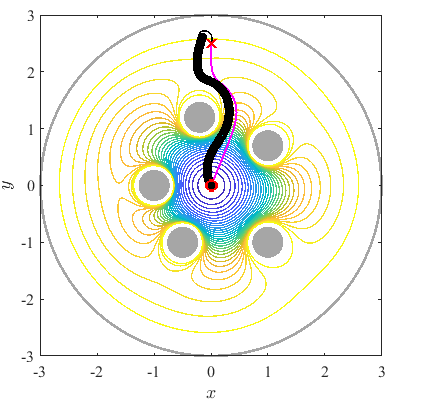}\label{fig:SIM_01_xy}}
	\caption{\label{fig:SIM_01}  Static source with unknown position. The level sets of the NF are plotted considering that the source position is known. The particle's trajectory estimated by the extremum seeking is in black color and the true gradient of the NF is in magenta color; the "x" marker is the starting point $(0, 2.5)$ and the read "o" is the source location $(0, 0)$.}
\end{figure}

We first present the simulation considering that the source is static and located at position $(0,0)$, as shown in Fig. \ref{fig:SIM_01}. For the sake of presentation, the level sets of the exact NF and its true gradient are also plotted.

As can be seen, the source was found by the robot and the resulting path was very close to the true gradient. 

\begin{figure}[hbt!]
	\centering
	\subfigure[][]{\includegraphics[width=0.7\columnwidth]{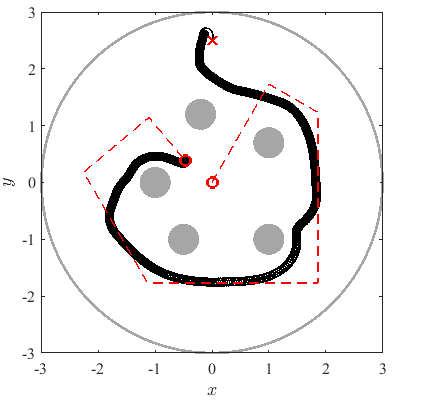}\label{fig:SIM_02_xy_vel_0_2}}
	\vspace{-0.2 cm}

	\caption{\label{fig:SIM_02_vel_0_2}  Moving source with $v_l = 0.2$ [m/s]. The robot's trajectory is in black color, the "x" marker is the starting point and the read "o" stands for the initial and final source location. The robot finds the source and tracks it up to the final position. 
	}
\end{figure}

Fig. \ref{fig:SIM_02_vel_0_2} considers the case of a slow time varying source. The source starts static and, after for 10 s, it follows the dashed red line path with speed $v_l = 0.2$ [m/s]. The source stops at position $(-0.47, 0.38)$ and keeps there until the end of the simulation. The particle finds the source while it is moving and tracks it up to the final position.

\subsection{Quadruple Mecanum Robot}

\begin{figure}[hbt!]
	\centering
	\includegraphics[width=0.4\linewidth]{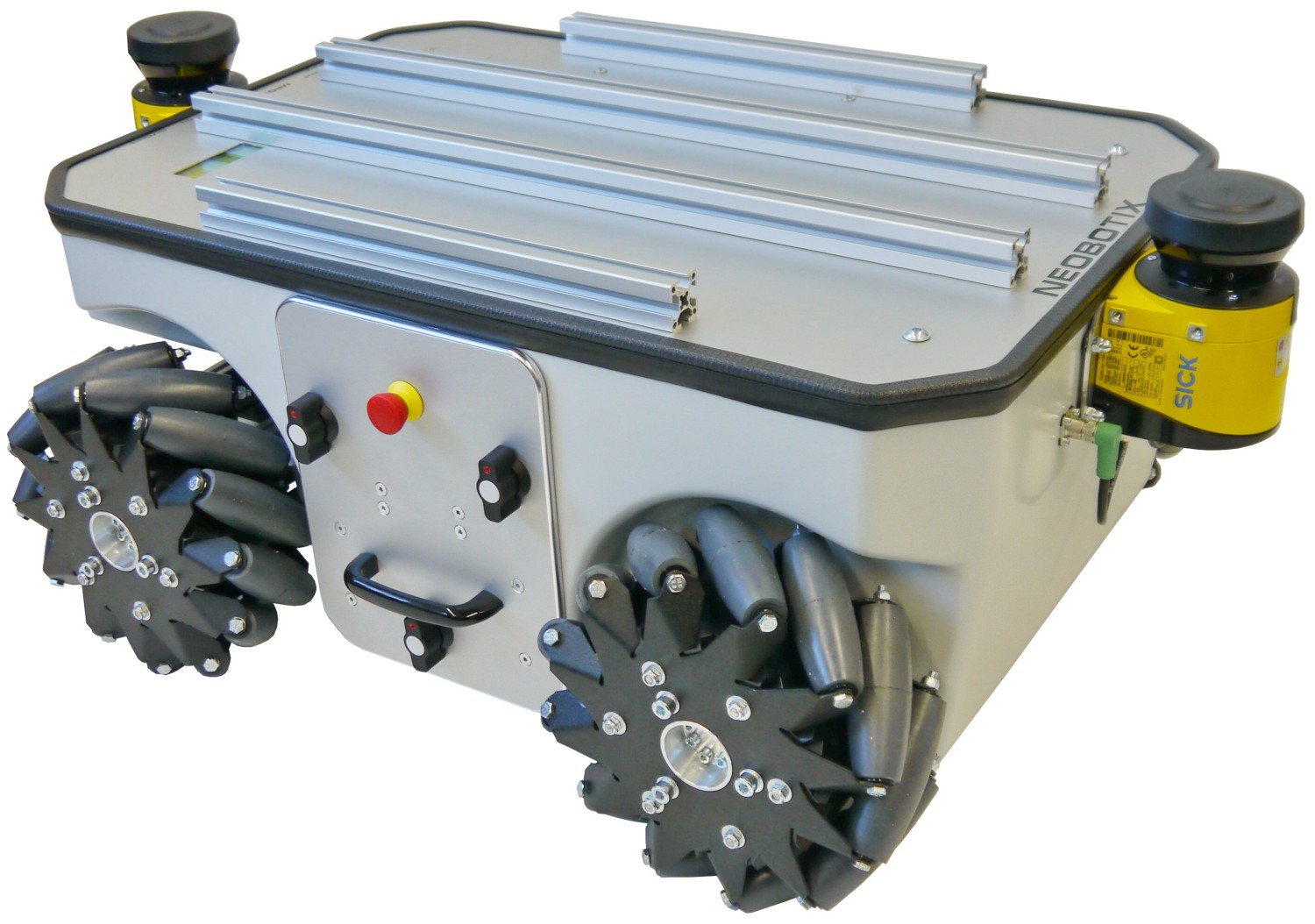}
	\caption{Robot Neobotix, model MPO-500 [source: neobotix-robots.com].}	
	\label{fig:MPO_500}
\end{figure}

A more realistic simulation was also carried out considering the Mecanum robot MPO-500, from Neobotix, which is a commercial robot with four mecanum wheel, as depicted in Fig. \ref{fig:MPO_500}. It is ROS compatible and has a Gazebo model. The main features of this robot are presented in Table \ref{tab:tab_param_mecanum}.

\begin{table}
	\caption{\label{tab:tab_param_mecanum} Parameters of the Neobotix Mecanum robot MPO-500.}
	\begin{center}
		\begin{tabular}{ cc } 
			\hline
			Parameter & value \\ 
			\hline
			Payload Default & 50 [kg]  \\ 
			Dimensions ($L\times W\times H$)  & 986 $\times$ 662 $\times$ 409 [mm] \\ 
			Max. linear speeds ($v_x$ and $v_y$) & 0.8 [m/s] \\ 
			Uptime & up to 7 h or up to 3 km\\
			Sensors & 1 or 2  2D laser scanner\\
			&  8 x ultrasonic sensors (optional)\\
			\hline
		\end{tabular}
	\end{center}
\end{table}

Simulations were performed in ROS+Gazebo using Rospy. ESC parameters were chosen such that $\omega = 2.5\pi$ rad/s, $\alpha = 0.15$, $C_x = C_y = -200$ and $h=\omega/2$. These values violate the gain assumption in Theorem \ref{theo:theo_Luiz}. The world consists in a circled environment with radius 7 m and five round obstacles, each one represented as $(x,\, y,\, r)$, being $r$ the radius. The five obstacles are $(4.5, \, 6.0, \, 0.35)$, $(6.5, \, 3.0, \, 0.7)$, $(6.0, \, 8.5, \,0.7)$, $(10.5, \, 7.5, \, 0.35)$ and $(11.0, \, 4.0, \, 0.35)$. The starting point of the robot is at $(2.0, 6.0)$.

Three simulations were considered. In the first one, all the obstacles are previously known, while in the second one, three out of five are known. In these two cases, it was considered $q_x =  q_y = 1.0$, and for the NF, it was set $k = 5$. 

In the third result, it is assumed that none of the obstacles is known \emph{a priori} and $k$ starts equal to zero. The initial potential (with $k=0$) is set as $0.1J/(0.1J + \beta_0)$. When a new obstacle is detected, the value of $k$ is incremented by one. 

The ESC algorithm was discretized with a sampling frequency equal to  $\omega_s = 10\omega$.  The maximum linear speed was limited to 0.8 m/s.
A rosbag file was generated during the simulations to record the ROS topics.

\begin{figure}[hbt!]
	\centering
	\includegraphics[width=0.7\linewidth]{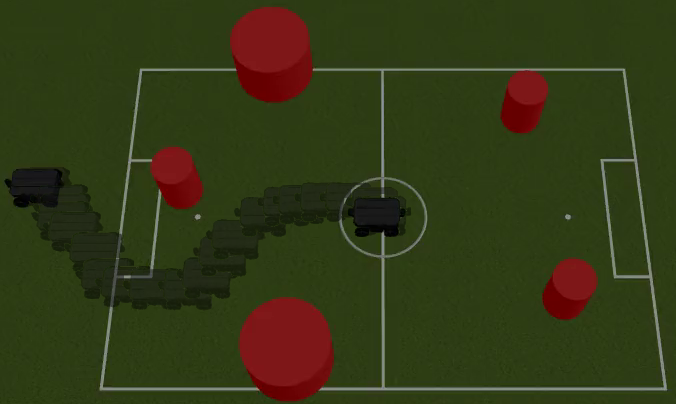}
	\caption{Gazebo simulation considering the Neobotix MPO-500. It is assumed that all the obstacles are previously known. 
		The source location is at the mid-field.}	
	\label{fig:GAZEBO_KNOWALL}
\end{figure}

\begin{figure}[hbt!]
	\centering
	\includegraphics[width=0.7\linewidth]{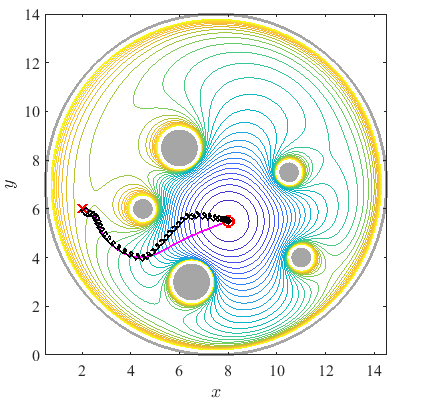}
	\caption{Plot of the robot's trajectory. The level sets of the NF are plotted considering that the source position is known. 
	}	
	\label{fig:KNOWALL_LS}
\end{figure}

Figs. \ref{fig:GAZEBO_KNOWALL} and \ref{fig:KNOWALL_LS} present the results when all the obstacles are known. It is possible to see that the trajectory of the robot follows the true gradient of the navigation function and the robot finds the source without hitting any obstacles.

\begin{figure}[hbt!]
	\centering
	\includegraphics[width=0.7\linewidth]{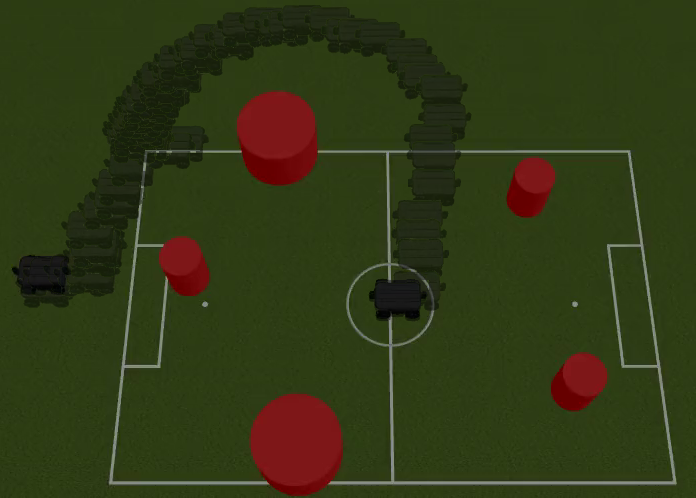}
	\caption{Gazebo simulation considering the Neobotix MPO-500. Three out of five obstacles are known.
	}	
	\label{fig:GAZEBO_KNOWPART}
\end{figure}

\begin{figure}[hbt!]
	\centering
	\includegraphics[width=.75\linewidth]{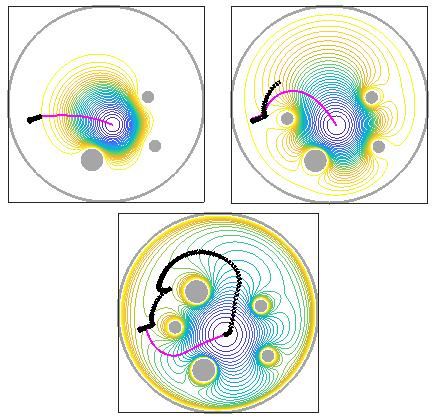}
	\caption{Trajectory of the robot on the level sets of the NF in three different moments, right before the updates of the NF and the final trajectory.
	}	
	\label{fig:KNOWPART_LS}
\end{figure}

Next results consider that the obstacles $(6.5, \, 3.0, \, 0.7)$, $(10.5, \, 7.5, \, 0.35)$ and $(11.0, \, 4.0, \, 0.35)$ are known \emph{a priori}. Fig. \ref{fig:GAZEBO_KNOWPART} shows the robot's trail in the Gazebo world, while Fig.  \ref{fig:KNOWPART_LS} presents the trajectory of the robot on the level sets of the NF, in three different moments: two right before the navigation function being updated, and the final trajectory. When a new obstacle is found, the NF is updated. The robot successfully finds the source, avoiding all the obstacles in the way. 

\begin{figure}[hbt!]
	\centering
	\includegraphics[width=0.7\linewidth]{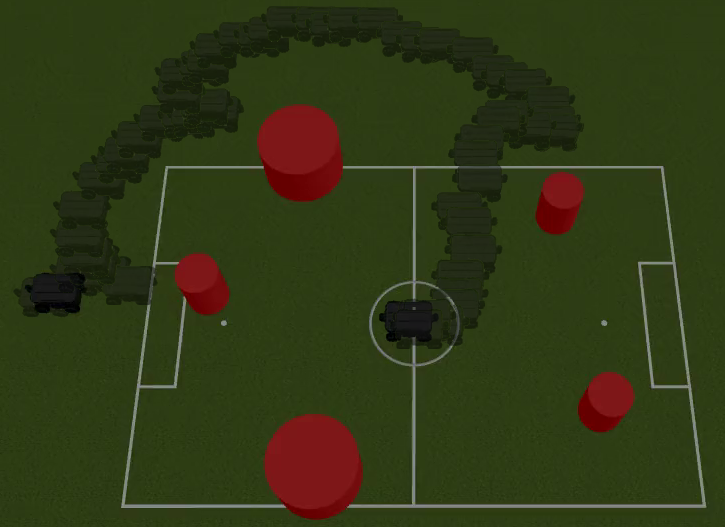}
	\caption{Gazebo simulation with Neobotix MPO-500, assuming that none obstacles are known \emph{a priori}.}	
	\label{fig:GAZEBO_KNOWNOT}
\end{figure}

\begin{figure}[hbt!]
	\centering
	\includegraphics[width=.9\linewidth]{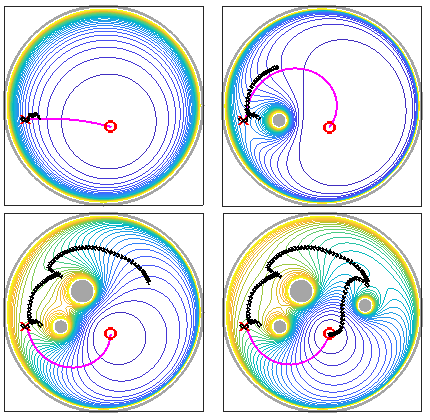}
	\caption{Trajectory of the robot on the level sets of the NF in four different moments, being three right before the updates of the NF, and the final trajectory. Obstacles are unknown \emph{a priori}. 
	}	
	\label{fig:KNOWNOT_LS}
\end{figure}

Next and more interesting result considers a completely unknown environment. Initially, $k=0$, but it increases by one when a new obstacle is found. Fig. \ref{fig:GAZEBO_KNOWNOT} shows the trail of the robot in the Gazebo, while Fig. \ref{fig:KNOWNOT_LS} presents the trajectory on the level sets of the NF. Again, the robot was able to find the source and avoid every obstacle in the way.

Even not satisfying the assumption $C_x, C_y \ll h$ in Theorem \ref{theo:theo_Luiz}, in all case, the robot was able to find the source and avoid all obstacles. We set high values for the loop gains because of the potential normalization imposed by the NF. Moreover, it is worth mentioning that the maximum speed of the robot was limited to 0.8 [m/s].

\subsection{Triple Omniwheel Robot}
	
\begin{figure}[hbt!]
	\centering
	\includegraphics[width=0.4\linewidth]{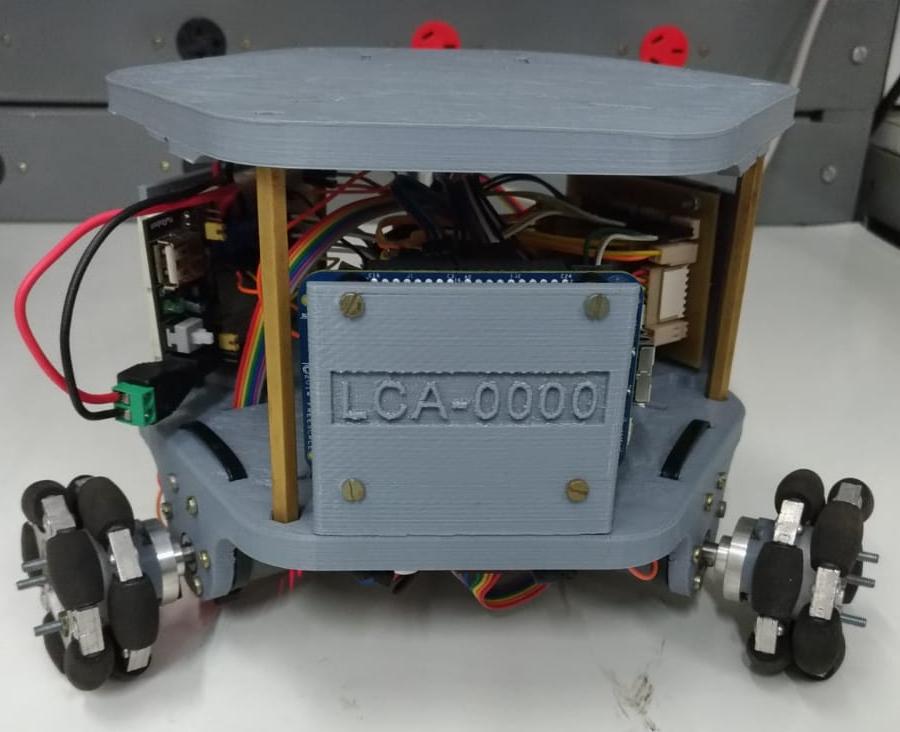}
	\caption{Schematic of the Triple Omniwheel Robot.}	
	\label{fig:OMNI}
\end{figure}

The parameters of the ESC were chosen such that $\omega = 2.5\pi$ rad/s, $\alpha = 0.15$, $C_x = C_y = -20$ and $h=\omega/2$. The world consists in a circled environment with radius 2.2 m and five round obstacles with radius 0.13 m at the points $(1.0, 0.0)$, $(1.5, 1.0)$, $(1.5, -1.0)$, $(2.5, 0.75)$, $(2.5, -0.75)$ , respectively. The experiment was performed with the real robot present in Figure \ref{fig:OMNI}, but the measurements of the source potential was artificially generated, according to Eq. (\ref{eq:eq_J}), using odometry. As in the last case of the previous subsection, it is assumed that the obstacle are unknown. However, since the robot has no lidar yet, there is a table with all the obstacles positions and it is assumed that when the robot gets close to an obstacle, its position and curvature are instantaneously detected and, thus, the NF is updated. That is the reason why we called this experiment as real dynamics with robot-in-the-loop. The ESC algorithm was discretized with a sampling frequency equal to  $\omega_s = 10\omega$.  The maximum $x$ and $y$ speeds were limited to 0.6 m/s. 

\begin{figure}[htbp]
	\center
	\subfigure[][]{\includegraphics[width=5.5cm]{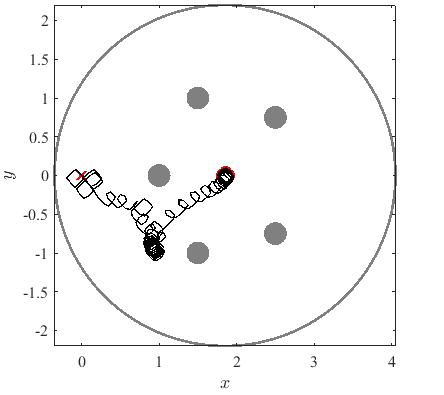}}
	\qquad
	\subfigure[][]{\includegraphics[width=5.5cm]{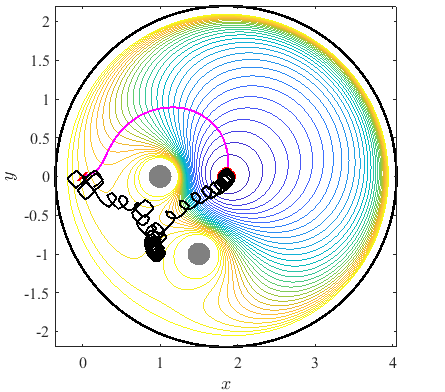}}
	\caption{\label{fig:Fig_pract} Results with robot in the loop. In a) it is presented the robot's trajectory (in black) considering the whole scenario, while b) shows the same trajectory with the obstacles detected by the robot and the level sets of final estimated navigation function with its true gradient (in magenta).}
\end{figure}

This robot was built for this experiment using a FRDM-K64F micro controller board from NXP. Its main parameters are shown in Table \ref{tab:tab_param_omni}. Practical result is presented in Figure \ref{fig:Fig_pract}. Notice that though it violates the assumptions of Theorem \ref{theo:theo_Luiz}, the robot could still reach the source collision-free.

\begin{table}
	\caption{\label{tab:tab_param_omni} Parameters of the triple omniwheel robot.}
	\begin{center}
		\begin{tabular}{ cc } 
			\hline
			Parameter & value \\ 
			\hline
			Wheels' radius ($R$) & 0.03 [m]  \\ 
			Distance from center to the wheels ($L$) & 0.15 [m] \\ 
			Robot's max. linear speeds ($v_x$ and $v_y$) & 0.6 [m] \\ 
			Robot's angular speed ($\Omega_R$) & 0.0 [rad/s]\\
			\hline
		\end{tabular}
	\end{center}
\end{table}

\section{Conclusions}\label{sec:sec_conclusions}

This paper addressed the problem of source seeking with scalar measurements in spaces cluttered with convex obstacles. To do so, the autonomous agent builds a navigation function online based on a local measurement of the source potential and its~(partial) knowledge of obstacles. A sinusoidal ESC loop was then used to minimize this artificial potential so as to smoothly navigate to the source location. Theoretical results guarantee that, under mild conditions, the robot attains its goal while avoiding collisions. Experiments with velocity-actuated particles and omni-directional robots corroborated the efficiency of the method, even in situations where no knowledge of the environment is available. Future work include testing this method with real potential measurements and extending it to 3D navigation for quadrotors or unmanned underwater vehicles.




\bibliographystyle{IEEEtran}
\bibliography{refs}

\section*{APPENDIX: Proof of Theorem~\ref{theo:theo_Luiz}}


\begin{proof}
Assuming that~$C_x = C_y = C$, start by noticing that the dynamical system in Fig.~\ref{fig:ES_NV} can be written as
\begin{subequations}\label{E:es_ode2}
\begin{align}
	\dot{\bp} &= -C \bz(\omega t) \big( \bar{\varphi}(\bp,\omega t) - \eta(t) \big)
		\label{E:es_ode_p2}
	\\
	\dot{\eta}(t) &= h \big[ \bar{\varphi}(\bp,\omega t) - \eta(t) \big]
		\label{E:es_ode_eta2}
\end{align}
\end{subequations}
where~$\bp = (x,y)^T$ denotes the particle position, $\bz(u) = \left[ \sin(u)\ -\cos(u) \right]^T$ is the modulation vector, and
\begin{equation}\label{E:varphibar}
	\bar{\varphi}(\bp,u) = \varphi\big( x + \alpha \sin(u), y - \alpha \cos(u) \big)
		\text{.}
\end{equation}
We know from~\cite[Thm.~1]{Zhang2007} that~\eqref{E:es_ode2} will converge to a neighborhood of~$\bp^\star$. Explicitly, $\limsup_{t \to \infty} \norm{\bp(t) - \bp^\star}^2 = \calO(\alpha^2 + 1/\omega^2)$. All that remain is therefore to show that~$\bp(t) \in \calF$ for all~$t \geq 0$. To do so, we proceed by finding a~$\calO(\alpha C/\omega)$ approximate description of the trajectory described by~\eqref{E:es_ode2} and show that it does not intersect any obstacle. Hence, \eqref{E:es_ode2} describes a collision-free trajectory as long as the agent never gets too close to the obstacles.

To obtain an approximate solution of~\eqref{E:es_ode2}, we start by using the classical multi-timescale approximation~\cite{O_Malley}. Since~$C \ll h$, there is an initial layer~$t = \calO(C/h)$ during which~$\bp$ is nearly constant and close to its initial value while~$\eta(t)$ converges exponentially fast to the quasi steady-state~$\eta_0(\bp(0))$ for
\begin{equation}\label{E:eta_closed2}
	\eta_0(\bp) = h e^{-ht} \int e^{ht} \bar{\varphi}(\bp,\omega t) dt
		\text{.}
\end{equation}
Then, for~$t = \calO(1)$, the value of~$\eta$ is close to the slow manifold, i.e., $\eta(t) = \eta_0(\bp(t)) + \calO(C/h)$, and~$\bp(t)$ satisfies~\eqref{E:es_ode_p2} with~$\eta(t) = \eta_0(\bp(t))$. Explicitly,
\begin{equation}\label{E:ode_p_outer}
	\dot{\bp} = -C \bz(\omega t) \big( \bar{\varphi}(\bp,\omega t) - \eta_0(\bp) \big)
		\text{,}
\end{equation}
for~$t = \calO(1)$. Notice that~$\eta$ from~$\bp$ are now decoupled.

To proceed, we apply the averaging method to obtain an approximate solution to~\eqref{E:ode_p_outer}. Start by normalizing the timescale using~$\tau = \omega t$ to get
\begin{equation}\label{E:ode_p_tau}
	\omega \dot{\bp} = -C \bz(\tau) \big( \bar{\varphi}(\bp,\tau) - \eta_0(\bp) \big)
		\text{.}
\end{equation}
Then, using the fact that~$\varphi \in \calC^2$, we can expand~\eqref{E:varphibar} around~$\bp$ as in
\begin{equation}\label{E:phi_taylor}
	\bar{\varphi}\big( \bp, \tau \big) = 
		\varphi(\bp) + \alpha \nabla \varphi(\bp)^T \bz(\tau) +
		\frac{\alpha^2}{2} \bz(\tau)^T \nabla^2 \varphi(\tilde{\bp}) \bz(\tau)
			\text{,}
\end{equation}
for some~$\bp - \alpha \ones \preceq \tilde{\bp} \preceq \bp + \alpha \ones$, where~$\preceq$ indicates the elementwise partial ordering, $\ones$ is a vector of ones, and~$\nabla f$ and~$\nabla^2 f$ denote the gradient vector and Hessian matrix of the function~$f$, respectively. Using~\eqref{E:phi_taylor} in~\eqref{E:eta_closed2} and~\eqref{E:ode_p_tau}, we obtain an expression of the form~$\dot{\bp} = \varepsilon g(\bp,\tau)$, where~$g$ is $2\pi$-periodic in~$\tau$ and~$\varepsilon = \frac{\alpha C}{\omega}$. Hence, 
\begin{equation}\label{E:avg_dynamics}
	\bp(\tau) = \bp_{av}(\tau) + \calO(\varepsilon)
		\text{,}
\end{equation}
where~$\dot{\bp}_{av} = \varepsilon g_{av}(\bp_{av})$ for the average dynamical system~\cite[Thm.~10.4]{Khalil}
\begin{equation*}
	g_{av}(\bp) = \frac{1}{2\pi} \int_{0}^{2\pi} g(\bp,\tau) d\tau
		\text{.}
\end{equation*}
Explicitly,
\begin{equation}\label{E:es_ode_av2}
	\dot{\bp}_{av} = -\frac{\alpha \omega C}{2(h^2 + \omega^2)}
		\left[ \bI + \frac{h}{\omega} \bJ \right] \nabla \varphi(\bp_{av})
		\text{.}
\end{equation}
where~$\bI$ is the identity matrix and
\begin{equation*}
	\bJ =
	\begin{bmatrix}
		0 & 1 \\ -1 & 0
	\end{bmatrix}
		\text{.}
\end{equation*}

It is worth noting, that although the agent trajectory~\eqref{E:avg_dynamics} follows~$\bp_{av}$ arbitrarily closely, $\bp_{av}$ is not the gradient dynamics~\eqref{E:dynamics} required to guarantee collision-free navigation due to the skew-symmetric term~$\bJ$. Despite this perturbation, this~\eqref{E:es_ode_av2} also avoid obstacles. Indeed, notice that since the obstacles~$\calO_i$ do not intersect, the value of~$\beta(\bp^\prime) = \prod_i \beta_i(\bp^\prime)$ determines whether the agent hits an obstacle. Explicitly,
\begin{gather*}
	\beta(\bp^\prime) \geq 0 \Leftrightarrow \bp^\prime \in \calF
	\text{, \ }
	\beta(\bp^\prime) = 0 \Leftrightarrow \bp^\prime \in \cup_i \del\calO_i
	\text{,}
	\\
	\text{and } \beta(\bp^\prime) < 0 \Leftrightarrow \bp^\prime \notin \calF
		\text{.}
\end{gather*}
Hence, if~$\beta$ is increasing along~$\bp_{av}$ whenever~$\bp_{av} \in \cup_i \del\calO_i$, i.e., whenever~$\beta(\bp_{av}) = 0$, then the dynamics~\eqref{E:es_ode_av2} is collision-free. To show this is the case, notice that
\begin{align*}
	\dot{\beta}(\bp_{av}) \Bigr|_{\beta(\bp_{av}) = 0}
		= \nabla\beta(\bp_{av})^T \dot{\bp}_{av} \Bigr|_{\beta(\bp_{av}) = 0}
		\text{,}
\end{align*}
which using~\eqref{E:es_ode_av2} yields
\begin{equation}\label{E:dotbeta}
\begin{aligned}
	\dot{\beta}(\bp_{av})
		= -\frac{\alpha \omega C}{2(h^2 + \omega^2)} \nabla\beta(\bp_{av})^T
		\left[ \bI + \frac{h}{\omega} \bJ \right] \nabla \varphi(\bp_{av})
		\text{,}
\end{aligned}
\end{equation}
where we omitted the fact that this expression is evaluated at~$\beta(\bp_{av}) = 0$ for clarity.

To proceed, we obtain from~\eqref{eq:eq_NV} that
\begin{multline}\label{E:gradphi}
	\nabla \varphi(\bp^\prime) =
		\left( f_0(\bp^\prime)^k + \beta(\bp^\prime) \right)^{-1-\frac{1}{k}}
	\\
	\left( \beta(\bp^\prime) \nabla f_0(\bp^\prime)
		- \frac{f_0(\bp^\prime)\nabla \beta(\bp^\prime)}{k} \right)
		\text{.}
\end{multline}
Using the value~\eqref{E:gradphi} in~\eqref{E:dotbeta} when~$\beta(\bp_{av}) = 0$, we obtain
\begin{equation}\label{E:dotbeta2}
\begin{aligned}
	\dot{\beta}(\bp_{av})
		= \frac{\alpha \omega C f_0(\bp_{av})^{-k}}{2k(h^2 + \omega^2)}
			\nabla\beta(\bp_{av})^T \left[ \bI + \frac{h}{\omega} \bJ \right]
				\nabla \beta(\bp_{av})
		\text{.}
\end{aligned}
\end{equation}
Observe from~\eqref{E:dotbeta2} that the sign of~$\dot{\beta}$ is determined by the quadratic form~$\nabla\beta(\bp_{av})^T \left[ \bI + \frac{h}{\omega} \bJ \right] \nabla \beta(\bp_{av})$, since the scalar factor is strictly positive. Notice, however, that since~$\bJ$ is skew-symmetric, the quadratic form~$\nabla\beta(\bp_{av})^T \bJ \nabla \beta(\bp_{av})$ vanishes. Thus, \eqref{E:dotbeta2} reduces to
\begin{equation}\label{E:dotbeta3}
\begin{aligned}
	\dot{\beta}(\bp_{av})
		= \frac{\alpha \omega C f_0(\bp_{av})^{-k}}{2k(h^2 + \omega^2)}
			\norm{\nabla\beta(\bp_{av})}^2 > 0
		\text{.}
\end{aligned}
\end{equation}
where we used the fact that~$\nabla\beta(\bp_{av})^T \nabla \beta(\bp_{av}) = \norm{\nabla\beta(\bp_{av})}^2$.

To conclude, \eqref{E:es_ode2} converges to a neighborhood of~$\bp^\star$~\cite[Thm.~1]{Zhang2007}. Moreover, the trajectory defined by~\eqref{E:es_ode_av2} is collision-free~[see~\eqref{E:dotbeta3}] and the true trajectory~$\bp(t)$ in~\eqref{E:avg_dynamics} followed by the agent is~$\calO(\alpha C/\omega)$ close to~\eqref{E:es_ode_av2}. Hence, if the obstacles are virtually enlarged by a factor of~$\calO(\alpha C/\omega)$, the agent will navigate to a neighborhood of~$\bp^\star$ while avoiding the actual obstacles.
\end{proof}




\end{document}